\documentclass[12pt]{article}
\usepackage{amsthm}

\begin{document}
\thispagestyle{empty}
\title{Proof of Aharoni Berger Conjecture} 
\date{}
\author{Vladimir Blinovsky}
\date{\small
 Institute for Information Transmission Problems, \\
 B. Karetnyi 19, Moscow, Russia,\\
vblinovs@yandex.ru}

\maketitle\bigskip

\begin{center}
{\bf Abstract}
We prove  Aharoni Berger Conjecture
\end{center}
\bigskip

{\bf Aharoni Berger Conjecture~(\cite{1})}
\bigskip

Let $B$ be properly $n$- colored bipartite multigraph with $n+1$ edges of each color. Then it contains rainbow matching of size $n$.
\bigskip

{\bf Proof}
\bigskip

Let $B=\{ V_1 =(v_{1,1},\ldots ,v_{1,p}),V_2 =(v_{2,1},\ldots ,v_{2,q}),E\}$ be bipartite graph as stated in Conjecture. W.l.o.g.
we assume that vertex $v_{11}$ is such that $|\{ K(E(v_{1,1},\cdot )\}|\neq n$. , where $K(E)$ is the color of edge $E$.

We define the following shifting procedure. Having graph $B$ we obtain new graph $$\bar{B}=\{\{ v_{1,1},\ldots ,v_{1,p}\},\{v_{2,1},\ldots, v_{2,q}\} ,\bar{E}\}$$ such that all edges in $B$ and $\bar{B}$ are the same with the following exceptions:

we have new edge $\bar{E}(v_{1,1},v_{2,\beta})$ instead of $E(v_{1p},v_{2,\beta})$ with the same color if no one edge $E(v_{1,1},\cdot)$ has this color and we have new edge $\bar{E}(v_{1,1},v_{2,\alpha })$ along with new edge $\bar{E}(v_{1,p},v_{2,\gamma})$ of the same color $K$ if edges $E(v_{1,1},v_{2,\gamma})$ along with  edge $E(v_{1,p},v_{2,\alpha})$ of the same color $K$ belong to $B$. 
It is easy to see that if graph $B$ does not contain rainbow matching of size  $n$ than $\bar{B}$ also does not have.

Continuing this shifting process we come to the bipartite $n$- colored multigraph $(V_1 ,V_2 ,E)$ which has $n+1$ edge of each color and $|V_1| =|V_2| =n+1$. 

Then we produce by induction: deleting some color $K$ and some vertex $v_{1,1}$, we obtain the subgraph, which using the same shifting procedure as before we reduce to the bipartite multigraph $(V_1 \setminus v_{1,1},V_2 \setminus v,E)$, where $v$ is the vertex in $V_2$ which has color $K$ and is connecting with vertex $v_{1,1}$. We obtain subgraph which has $n$ edge each color and $n$ vertices in each part. By induction it contais $n-1$-rainbow matching. Adding edge $(v,v_{1,1})$ we obtain 
$n$- rainbow matching of initial graph which is sufficient to our proof. Note that we start the  induction from $n=2$ in which case the statement of conjecture is obviously true.

This completes the proof.

\end{document}